\def\[{[\![}
\def\]{]\!]}
\def\l{\ldots}
\def\q{\quad}
\def\h{\hat{h}}
\def\e{\hat{e}}
\def\f{\hat{f}}
\def\n{\noindent}

\def\a{\hat{a}}
\def\k{{\bar k}}
\def\L{{\bar L}}
 \def\q{{\bar q}}
 \def\x{{\bar x}}

\hfill {\bf SISSA 123/96/FM}
\vskip 6mm

\baselineskip 24pt   

\noindent
\centerline{\bf  AN ALTERNATIVE TO THE CHEVALLEY DESCRIPTION OF 
U[sl(n+1)] AND U$_q$[sl(n+1)] }

\vskip 32pt
\noindent
\hskip 55mm Tchavdar D. Palev\footnote*{Permanent Address: 
 Institute for Nuclear Research and Nuclear Energy, 1784 Sofia, 
Bulgaria; $e$-mail: tpalev@inrne.acad.bg} and Preeti Parashar

\noindent
\centerline{International School for Advanced Studies, via Beirut 2-
4,
34013 Trieste, Italy} 
\smallskip

\vskip 48pt
\noindent
{\bf Abstract.} An alternative to the Chevalley description of
the universal enveloping algebra $U[sl(n+1)]$ and of its
$q-$deformed analogue $U_q[sl(n+1)]$ in terms of generators and
relations is given. In particular $U[sl(n+1)]$ is an associative
algebra with 1, generators $\a_1^\pm, \l,\a_n^\pm$ and
relations 
$[\a_1^\xi,\a_2^\xi]=0,\;
[[\a_i^\xi,\a_{j}^{-\xi}],\a_k^\xi]=\delta_{jk}\a_i^\xi+
  \delta_{ij}\a_k^\xi,
\; \break
|i-j|\le 1,  \; \xi=\pm .
$

\vskip 48pt

\vskip 48pt

\n {\bf 1. Introduction}

\bigskip
\n In the present paper we describe firstly the special linear
Lie algebra $sl(n+1)$ and its universal enveloping algebra
($UEA$) $U[sl(n+1)]$ via creation and annihilation generators
(CAGs) [1]. Thus we provide an alternative to the Chevalley
description of $U[sl(n+1)]$ in terms of generators and relations.
Secondly, we describe the quantized universal enveloping algebra
$U_q[sl(n+1)]$ [2, 3] entirely via deformed CAGs.

The concept of creation and annihilation operators (or
generators) of a (semi)simple Lie (super)algebra of rank $n$ was
introduced in [1]. The annihilation generators $\a_1^-,\l,\a_n^-$
(resp. the creation generators $\a_1^+,\l,\a_n^+$) are always
among the positive (resp. among the negative) root vectors. 
To be more precise we recall that the root system of
$sl(n+1)\equiv A_n$ is
$$
\Delta=\{ \varepsilon_A - \varepsilon_B
|A\neq B=0,1,\ldots,n \}.
\eqno(1)
$$
The correspondence between the CAGs of $sl(n+1)$ and their roots
reads: $\a_i^\pm \leftrightarrow \mp( \varepsilon_0 -
\varepsilon_i)$. The simple root vectors $\e_i$ have roots
$\varepsilon_{i-1} - \varepsilon_i$, $\e_i
\leftrightarrow \varepsilon_{i-1} - \varepsilon_i,\; i=1,\l,n $.
Therefore the CAGs are very different from the Chevalley
generators $\e_i,\f_i,\h_i,\;i=1,\l,n$ of $sl(n+1)$.

For us personally the interest in the subject stems from the
observation that the CAGs (including $sl(\infty)$) provide an
alternative way for quantization of integer spin fields.  This
leads to a generalization of the quantum statistics, the
$A$-statistics [4].  The latter, as it is clear now, is a
particular case of the Haldane exclusion statistics [5], a
subject of considerable interest in condensed matter physics
(see, for instance, [6-8] and the references therein).  We are
not going to discuss the properties of the underlying statistics
here. We note only that the quantization of $U[sl(n+1)]$ leads to
deformations of the Fock space representations, considered in
[4], which gives 
new solutions for the microscopic statistics
of the $g$-ons of Karabali and Nair [9], a particular realization
of the Haldane statistics.  The very fact that we consider Hopf
algebra deformations of the creation and the annihilation
operators has an additional advantage: using the comultiplication
one may construct new representations of $\a_1^\pm,\a_2^\pm,\l$
in any tensor product of state spaces, namely new solutions for
the $g$-on statistics.


\vskip 32pt
\noindent
{\bf 2. The Lie algebra $sl(n+1)$}
\bigskip

\noindent
In order to define the CAGs of $sl(n+1)$ it is convenient to
consider it as a subalgebra of the Lie algebra $gl(n+1)$.
The $UEA$ of $gl(n+1)$ can
be defined as an associative algebra with 1 of the indeterminates
$ \{e_{AB}|A,B=0,1,\ldots,n\}$ subject to the relations
(below and throughout $[x,y]=xy-yx$)
$$
[e_{AB},e_{CD}]=\delta_{BC}e_{AD}-\delta_{AD}e_{CB}.\eqno(2)
$$
Then $gl(n+1)$ is a subalgebra of the Lie algebra  $U[gl(n+1)]$
with generators $e_{AB},\;A,B=0,1,\ldots,n$ and commutation 
relations (2). 

The Cartan subalgebra $H'$ of $gl(n+1)$ has a basis 
$e_{00},e_{11},\l,e_{nn}$. Let $\varepsilon_0,\varepsilon_1,\l,
\varepsilon_n$ be the dual basis, $\varepsilon_A(e_{BB})=\delta_{AB}$.
The root vectors of both $gl(n+1)$ and $sl(n+1)$ are
$e_{AB},\;A\ne B=0,1,\l,n$.  The root of each $e_{AB}$ is 
$\varepsilon_A -\varepsilon_B$. Then
$$
sl(n+1)=span\{e_{AB}, e_{AA}-e_{BB}|A\ne B=0,1,\l,n\}. \eqno(3)
$$

The canonical description of $sl(n+1)$, which is appropriate for 
Hopf algebra deformations, is given in terms of its
Chevalley generators 
$$
\e_i=e_{i-1,i},\quad  \f_i=e_{i,i-1},\quad 
\h_i=e_{i-1,i-1}-e_{ii},\quad i=1,2,\l,n,
\eqno(4)
$$
and the $n\times n$ Cartan matrix $\{\alpha_{ij}\}$ with entries
$\alpha_{ij}=2\delta_{ij}-\delta_{i,j-1}-\delta_{i-1,j}$:
$U[sl(n+1)]$ is an associative algebra with 1 of the Chevalley
generators subject to the Cartan relations 
$$
[\h_i,\h_j]=0,\;
[\h_i,\e_j]=\alpha_{ij}\e_j,\;
[\h_i,\f_j]=-\alpha_{ij}\f_j,\;
[\e_i,\f_j]=\delta_{ij}\h_i, \eqno(5)
$$
and the Serre relations
$$
\eqalign{
& [\e_i,\e_j]=0,\quad [\f_i,\f_j]=0, \quad \vert i-j \vert>1, \cr
& [\e_i,[\e_i,\e_{i\pm 1}]]=0,\quad  
  [\f_i,[\f_i,\f_{i\pm 1}]]=0.\cr
}\eqno(6)
$$

The creation and the annihilation generators of $sl(n+1)$ are [4] 
$\a_i^+=e_{i0}$ and $\a_i^-=e_{0i},\; i=1,\l,n$. 
Our new result is contained in the following proposition.

{\it Proposition 1.} $U[sl(n+1)]$ is an associative algebra 
with 1, (free) generators $\a_1^\pm, \a_2^\pm,\l,\a_n^\pm$ and 
relations
$$
[\a_1^\xi,\a_2^\xi]=0,\;\;
[[\a_i^\xi,\a_{j}^{-\xi}],\a_k^\xi]=\delta_{jk}\a_i^\xi+
  \delta_{ij}\a_k^\xi, \quad 
|i-j|\le 1,  \;\; \xi=\pm,\;\; i,j,k =1,\l,n .\eqno(7)
$$

We skip the proof, since it is a particular case of the proof of the
Theorem in Sect.3 ,
when $q \rightarrow 1$. The expressions for the
Chevalley generators in terms of the CAGs read ($i\ne 1$):
$$
\eqalign{
& \e_1=\a_1^-,\quad \e_i=[\a_{i-1}^+,\a_i^-],   \cr
& \f_1=\a_1^+,\quad \f_i=[\a_{i}^+,\a_{i-1}^-], \cr
& \h_1=[\a_1^-,\a_1^+],\quad 
  \h_i=[\a_i^-,\a_i^+]-[\a_{i-1}^-,\a_{i-1}^+],\cr
}\eqno(8)
$$
with inverse relations
$$
\eqalign{
& \a_1^-=\e_1,\quad \a_i^-=[[[\l[[\e_1,\e_2],\e_3]\l],\e_{i-
1}],\e_i],\cr
& \a_1^+=\f_1,\quad \a_i^+=[\f_i,[\f_{i-1},[\l 
[\f_3,[\f_2,\f_1]]\l]]].\cr
}\eqno(9)
$$

\n Only from (5), (6) and (9) one derives all relations among the
CAGs:
$$
[[\a_i^\xi,\a_j^{-\xi}],\a_k^\xi]=\delta_{jk}\a_i^\xi+
  \delta_{ij}\a_k^\xi,\quad 
[\a_i^\xi,\a_j^\xi]=0, \quad \xi=\pm, \quad i,j,=1,\l,n. \eqno(10)
$$
Note that (7) is a subset of (10). This is to say that (7)
represents the minimal set of relations required to define
$U[sl(n+1)]$. The fact that the creation (annihilation)
generators commute among themselves facilitates the construction
of the Fock representations of $sl(n+1)$ (also in the case
$n=\infty$) [4]. A similar problem for the parastatistics of
an arbitrary order [10, 11] remains still unsolved: explicit
expressions for the transformations of the Fock spaces do not
exist.


\vskip 32pt
\noindent
{\bf 3. $U_q[sl(n+1)]$ in terms of deformed creation and 
        annihilation generators}
\bigskip

\noindent
The description of $U_q[sl(n+1)]$ in terms of Chevalley
generators is well known. See, for instance [2, 3], where all
Hopf algebra operations are explicitly written. Here we write only
the algebra operations. $U_q[sl(n+1)]$ is an associative
algebra with 1 of its Chevalley generators
$e_i,\;f_i,\; k_i=q^{\h_i},\;\bar{k_i}\equiv k_i^{-1},\;i=1,\l,n,$
subject to the Cartan relations ($\bar{q}\equiv q^{-1}$) 
$$
\eqalignno{
&  k_i\k_i=\k_ik_i=1, \quad  k_ik_j=k_jk_i, & (11a) \cr
&  k_ie_j=q^{\alpha_{ij}}e_jk_i,\quad k_if_j=q^{-\alpha_{ij}}f_jk_i,&
   (11b)\cr
&  [e_i,f_j]=\delta_{ij}{{k_i-{\bar k}_i}\over{q-{\bar q}}},& (11c) 
}
$$
and the Serre relations 
$$
\eqalignno{
& [e_i,e_j]=0, \quad  [f_i,f_j]=0,\quad \vert i-j \vert \neq 1,
  & (12a) \cr
& [e_i,[e_i,e_{i \pm 1}]_{\bar q}]_q=
[e_i,[e_i,e_{i \pm 1}]_q]_{\bar q}=0, \quad 
[f_i,[f_i,f_{i \pm 1}]_{\bar q}]_q=
[f_i,[f_i,f_{i \pm 1}]_q]_{\bar q}=0. & (12b)\cr 
}
$$
Here and throughout we assume that $q$ is not a root of 1.

Having in mind the expressions (9) we define the deformed
CAGs as follows (below and throughout $[a,b]_x=ab-xba$):
$$
\eqalignno{
& a_1^-=e_1,\quad a_i^-
  =[[[\l[[e_1,e_2]_{\q},e_3]_{\q}\l]_{\q},e_{i-1}]_{\q},e_i]_{\q}
  =[a_{i-1}^-,e_i]_{\q}, \quad i\ne 1, & (13a)\cr
& a_1^+=f_1,\quad a_i^+=[f_i,[f_{i-1},[\l 
[f_3,[f_2,f_1]_q]_q\l]_q]_q]_q
       =[f_i,a_{i-1}^+]_q,
 \quad i\ne 1. & (13b)\cr
}
$$

Note that Eqs. (11), (12) are invariant with respect to the
linear antiinvolution $(\;)^*$, defined as
$$
(e_i)^*=f_i, \; (k_i)^*=\k_i, \; (q)^*=\q , \; ((x)^*)^*=x, \; 
(ab)^*=(b)^*(a)^* .\eqno(14a)
$$
Therefore, if $F, G\in U_q[sl(n+1)]$ and $F=G$, then also 
$(F)^*=(G)^*$. In particular,
$$
(a_i^\pm)^*=a_i^\mp. \eqno(14b)
$$

{\it Proposition 2}. The following "mixed" relations hold ($i\ne 1$):
$$ 
\eqalign{
& (a) \;\;[e_i,a_j^-]_{q^{\delta_{i-1,j}-\delta_{ij}}}
      =-q\delta_{i-1,j}a_i^-,\quad\quad (b)\;\;
  [f_i,a_j^+]_{q^{\delta_{i-1,j}-\delta_{ij}}}=\delta_{i-1,j}a_i^+,\cr
& (c) \;\; [e_i,a_j^+]=\delta_{ij}a_{i-1}^+\k_i,\quad\quad
    \hskip 18mm (d) \;\;  [f_i,a_j^-]=-\delta_{ij}k_ia_{i-1}^-
,\cr}\eqno(15)
$$

{\it Proof.} The proof is based on repeated use of identities 
like ($\x=x^{-1}$)
$$
\eqalignno{ 
& If\;\; [a,b]=0,\;\; then \;\; [[a,c]_q,b]_p=[a,[c,b]_p]_q;& (16a)\cr
& If\;\; [a,c]=0,\;\;
then \;\; (x+\x)[b,[a,[b,c]_x]_x]=[a,[b,[b,c]_x]_{\x}]_{x^2}-
[[b,[b,a]_x]_{\x},c]_{x^2}. & (16b) \cr
}
$$
We begin with (15a).

\n (i) Let $j<i-1$. Then from $(12a)$ one immediately has
$
[e_i,a_j^-]=0. 
$

\n (ii) Let $j=i-1$. From (13$a$) $ [a_{i-1}^-,e_i]_{\q}=a_i^-$ and 
therefore
$
[e_i,a_{i-1}^-]_q=-qa_i^-. 
$

\n (iii) Let $j=i$. If $i=2$, $[e_2,a_2^-]_\q=[e_2,[e_1,e_2]_\q]_\q=
   -\q[e_2,[e_2,e_1]_q]_\q=0$ according to (12b). If $i>2$ set
$a_i^-=[[a_{i-2}^-,e_{i-1}]_\q, e_i]_\q$. Then
$[e_i,a_i^-]_\q=[e_i,[[a_{i-2}^-,e_{i-1}]_\q, e_i]_\q,e_i]_\q $ and, 
since $[e_i,a_{i-2}^-]=0$, using twice $(16a),$ one has:
$[e_i,a_i^-]_\q=[e_i,[a_{i-2}^-,[e_{i-1},e_i]_\q]_\q]_\q=
[a_{i-2}^-,[e_i,[e_{i-1},e_i]_\q]_\q]_\q=0 $, since   
$[e_i,[e_{i-1},e_i]_\q]_\q=-\q[e_i,[e_{i},e_{i-1}]_q]_\q=0$ 
according to (12b). Thus,
$
[e_i,a_i^-]_\q=0. 
$

\n (iv) Let $i+1=j$. From $(16)$ one has at $x=\q$  

\n
$$[e_i,[[e_{i-1},e_i]_\q,e_{i+1}]_\q]
=(q+\q)^{-1}([e_{i-1},[e_i,[e_i,e_{i+1}]_\q]_q]_{\q^2}
-[e_i,[e_i,e_{i-1}]_\q]_q,e_{i+1}]_{\q^2})=0,\eqno(17)
$$
according to (12b). 
Therefore $[e_2,a_3^-]=[e_2,[[e_1,e_2]_\q,e_3]_\q]=0 $.
If $i>2$, $a_{i+1}^-=[[[a_{i-2}^-,e_{i-1}]_\q,e_i]_\q,e_{i+1}]_\q  $
and since $a_{i-2}^-$ commutes with $e_i, e_{i+1}$, from (16a) one
obtains $a_{i+1}^-=[a_{i-2}^-,[[e_{i-1},e_i]_\q,e_{i+1}]_\q]_\q$.
Therefore,
$
[e_i,a_{i+1}^-]=[e_i,[a_{i-2}^-,[[e_{i-1},e_i]_\q,e_{i+1}]_\q]_\q ]=
 [a_{i-2}^-,[e_i,[[e_{i-1},e_i]_\q,e_{i+1}]_\q] ]_\q=0, 
$ 
according to (17).  

\n (v) Let $j>i+1$. Then
$a_j^-=[[\l[[a_{i+1}^-,e_{i+2}]_\q,e_{i+3}]_\q \l ]_\q,e_{j-
1}]_\q,e_j]_\q$
and, since $e_i$ commutes with $e_{i+2},\;e_{i+3},\l e_{j}$, according
to (16a), and it commutes also with $a_{i+1}^-$, according to (iv), 
one
has
$
[e_i,a_j^-]=0. 
$
The unification of (i)-(v) yields $(15a)$.   
Applying the antiinvolution (14) on both sides of 
$(15a)$ one obtains $(15b)$. 

\n We pass to prove $(15c)$.

\n (i) For $i>j$ , $(15c)$ is an immediate consequence of $(11c)$ and 
$(13b)$.

\n (ii) Let $i=j$. $[e_i,a_i^+]=[e_i,[f_i,a_{i-1}^+]_q]$ (from
(i) and $(16a)$) $=[[e_i,f_i],a_{i-1}^+]_q $
\hfill\break
$=[(q-\q)^{-1}(k_i-\k_i,)a_{i-1}^+]_q=a_{i-1}^+\k_i.$ 

\n (iii) Let $j=i+1$. Using $(16a)$ and then $(15a)$, we have

\n
$[e_i,a_{i+1}^+]=[e_i,[f_{i+1},[f_i,a_{i-1}^+]_q]_q]
= [f_{i+1},[[e_i,f_i],a_{i-1}^+]_q]_q $
\hfill\break 
$=[f_{i+1},[(q-\q)(k_i-\k_i),a_{i-1}^+]_q]_q=[f_{i+1},a_{i-
1}^+\k_i]_q = 0.
$

\n (iv) For $j>i+1$,
$
a_j^+=[f_{j},[f_{j-1},[\l [f_{i+3},[f_{i+2},a_{i+1}^+]_q]_q\l]_q]_q]_q
$
and since $e_i$ commutes with $a_{i+1}^+$ according to (iii) and
with $f_{j}, f_{j-1},\l, f_{i+2}$ 
according to $(11c)$, we conclude that $[e_i,a_{j}^+]=0$. 

\n Combining (i)-(iv) gives $(15c)$. Acting on both sides
of it with the antiinvolution $(14)$ one derives $(15d)$.

{\it Proposition 3.}  The deformed CAGs $(13)$
together with the "Cartan" operators $k_1,\ldots,k_n$
generate (in a sense of an associative algebra)
$U_q[sl(n+1)]$.

{\it Proof.} The proof is an immediate consequence of the
relations:
$$
\eqalignno{
& [a_i^-,a_i^+]={{k_1k_{2}\ldots k_i -
{\bar k}_1{\bar k}_{2}\ldots {\bar k}_i}\over q-{\bar q}},\quad
i=1,\ldots,n,& (18a) \cr 
& [a_i^-,a_{i+1}^+]=-k_1k_{2}\ldots k_{i}f_{i+1}, \quad
[a_{i+1}^-,a_i^+]=-e_{i+1}{\bar k}_1{\bar k}_{2}\ldots {\bar
k}_{i}, \quad i=1,\ldots,n-1.& (18b) \cr 
}
$$
These equations are proved by induction on $i$. For $i=1$ $(18a)$ 
holds.
Let $(18a)$ be true. Then from (13$b$), $(16a)$ and $(18a)$

\n $[a_i^-,a_{i+1}^+]=[a_i^-,[f_{i+1},a_i^+]_q]=
[f_{i+1},[a_i^-,a_i^+]]_q={1\over{q-\q}}
[f_{i+1},k_1k_2\l k_i-\k_1\k_2\l \k_i]_q $.

\n Using (11) the latter yields
 
\n $[a_i^-,a_{i+1}^+]=-qf_{i+1}k_1k_2\l k_i
 =-k_1k_2\l k_if_{i+1}$.

\n Similarly one derives the other equation in $(18b)$. The 
conclusion 
so far is that if $(18a)$ holds, then also Eqs. $(18b)$ hold. 
Assuming 
this, compute
$[a_{i+1}^-,a_{i+1}^+]=[[a_i^-,e_{i+1}]_\q,a_{i+1}^+]$. Using the 
identity  $[[A,B]_x,C]=[[A,C],B]_x + [A,[B,C]]_x$, one has

\n
$[a_{i+1}^-,a_{i+1}^+]=[[a_{i}^-,a_{i+1}^+],e_{i+1}]_\q +
[a_i^-,[e_{i+1},a_{i+1}^+]]_\q $
\hfill\break
$=-[k_1k_2\l k_if_{i+1},e_{i+1}]_\q
+[a_i^-,a_i^+\k_{i+1}]_\q $
\hfill\break
$=[e_{i+1},f_{i+1}]k_1k_2 \l k_i
+[a_i^-,a_i^+]\k_{i+1} $
\hfill\break
$=(q-\q)^{-1}(k_1k_2\l k_{i+1}-
\k_1\k_2 \l \k_{i+1})
$.

\n Thus, if $(18a)$ holds for a certain $i$, it holds for any $i$. 
Hence
Eqs. $(18a)$ and $(18b)$ hold. Therefore $(i=1,2,\l ,n-1)$,
$$
e_1=a_1^-,\quad f_1=a_1^+,\quad
e_{i+1}=-[a_{i+1}^-,a_i^+]k_1k_2\l k_i,\quad
f_{i+1}=-\k_1\k_2\l \k_i[a_i^-,a_{i+1}^+].
\eqno(19)
$$
Hence $a_1^\pm,\l ,a_n^\pm $ together with $k_1^{\pm
1},\l,k_n^{\pm 1}$ generate $U_q[sl(n+1)]$, which completes the
proof.

Define new ``Cartan'' generators
$$
L_i=k_1k_2\l k_i \Leftrightarrow k_1=L_1,\;\; k_i=L_iL_{i-1}^{-1},
\quad i\ne 1. \eqno(20)
$$
We are now ready to formulate our main result.

{\it Theorem.} $U_q[sl(n+1)]$ is an associative algebra with unity,
generators $a_i^\pm,\;L_i,\;L_i^{-1}\equiv {\bar L}_i,\;i=1, \l ,n$
and relations
$$
\eqalignno{
&  L_i\L_i=\L_iL_i=1, \quad  L_iL_j=L_jL_i, & (21a)\cr
& L_ia_j^\pm=q^{\mp (1+ \delta_{ij})}a_j^\pm L_i, & (21b)  \cr
& [a_i^-,a_i^+]={L_i-{\bar L}_i\over q-{\bar q}},& (21c) \cr 
& [[a_i^{\eta},a_{i+ \xi}^{-\eta}],  
a_j^{\eta}]_{q^{\xi(1+\delta_{ij})}}
=\delta _{j,i + \xi}L_j^{-\xi \eta}a_i^{\eta}
 \quad \xi , \; \eta=\pm \;or \; \pm1,& (21d) \cr 
& [a_1^\eta,a_2^\eta]_q=0,\quad \eta=\pm  .& (21e) \cr 
}
$$

{\it Proof.} Most of the preliminary results, necessary for the
proof, are already obtained.
Eq. $(21a)$ is an immediate consequence of $(11a)$ and the definition
(20). From (9), (11$b$) and (20) one derives 
$
L_ia_j^+=q^{\sum_{r=1}^i\sum_{s=1}^j \alpha_{rs}}a_j^+L_i.
$
Then (21b) follows from the observation that
$
\sum_{r=1}^i\sum_{s=1}^j\alpha_{rs}=1+\delta_{ij}. 
$
Inserting (20) in $(18a)$ one obtains $(21c)$. Replacing in (15) 
$e_i$ and
$f_i$ with the right hand sides of (19) one derives all triple
relations $(21d)$. The nontrivial part is to write all of them in
the compact form $(21d)$. Eqs. $(21e)$ coincide with two of the
Serre relations $(12b)$.

Thus, Eqs. (21) hold. It remains to prove that these relations
are sufficient in order to derive any other relation in
$U_q[sl(n+1)]$. To this end it suffices to show that the Cartan
relations (11) and the Serre relations (12) follow from Eqs.
(21). The relations for the Chevalley generators are 
(19) and (20). Eq. $(11a)$ follows straight from $(21a)$, and
$(11b)$ is obtained by the repeated use of $(21b)$. 

The proof of Eq. $(11c)$ is however not so simple. We consider in
detail the more difficult case, namely when $i,j\ne 1$. Then
\n
$[e_i,f_j] = [[a_{i-1}^+,a_i^-]L_{i-1},\L_{j-1}[a_j^+,a_{j-1}^-]] $
\hfill\break
$= [[a_{i-1}^+,a_i^-],\L_{j-1}][a_j^+,a_{j-1}^-]L_{i-1} + 
\L_{j-1}[[a_{i-1}^+,a_i^-],[a_j^+,a_{j-1}^-]]L_{i-1} $
\hfill\break
$+[a_{i-1}^+,a_i^-]\L_{j-1}[L_{i-1},[a_j^+,a_{j-1}^-]]$ since 
$[L_{i-1},\L_{j-1}] = 0.$ 

\n Let us shift all the $L$ and $\L$
to the left by making repeated use of Eq.$(21b)$. So r.h.s. becomes

\n
$ q^{1-\delta_{i-1,j-1}+\delta_{i,j-1}} \L_{j-1} L_{i-1}
([a_{i-1}^+,a_i^-][a_j^+,a_{j-1}^-] - q^{\delta_{i-1,j}-\delta_{i,j-
1}} 
[a_j^+,a_{j-1}^-][a_{i-1}^+,a_i^-])$, 
\n which can be writen as 
$$
[e_i,f_j] =q^{1-\delta_{i-1,j-1}+\delta_{i,j-1}} \L_{j-1}L_{i-1}
[[a_{i-1}^+,a_i^-],[a_j^+,a_{j-1}^-]]_{q^{\delta_{i-1,j}-\delta_{i,j-
1}}}.
\eqno(22)
$$
Now we use the identity 
$$ 
[A,[B,C]_x]_y = [[A,B]_z,C]_t + z[B,[A,C]_r]_s \quad  with \; 
conditions \;\; x = zs,\; y = zr,\; t = zsr. \eqno(23)
$$ 
\n (i) For $i = j$ (22) reduces to 
$[e_i,f_i] = [[a_{i-1}^+,a_i^-],[a_i^+,a_{i-1}^-]]=
\hfill\break
[[[a_{i-1}^+,a_i^-],a_i^+]_q, a_{i-1}^-]_\q - 
q[a_i^+,[[a_i^-,a_{i-1}^+],a_{i-1}^-]_\q]_\q  = 
[\L_ia_{i-1}^+,a_{i-1}^-]_\q - q[a_i^+,\L_{i-1}a_i^-]_\q $ (using 
$(21d)$)\hfill\break 
$= \L_i[a_{i-1}^+,a_{i-1}^-] - \L_{i-1}[a_i^+,a_i^-]$ (from 
$(21b)$) $ = {(L_i\L_{i-1} - L_{i-1}\L_i) \over q - {\bar q}}= 
{{k_i-{\bar k}_i}\over{q-{\bar q}}}$.

\n (ii) For $|i-j|>1$ (22) reduces to 
$[e_i,f_j]=q\L_{j-1}L_{i-1}[[a_{i-1}^+,a_{i}^-],[a_j^+,a_{j-1}^-]]
\hfill\break
=q\L_{j-1}L_{i-1}\{[[[a_{i-1}^+,a_{i}^-],a_j^+]_q,a_{j-1}^-]_\q
-q[a_j^+,[[a_{i}^-,a_{i-1}^+],a_{j-1}^-]_\q]_\q\}$ $=0$ according to 
$(21d)$.

\n (iii) For $j=i-1$ (22) reduces to
$[e_i,f_{i-1}] =q\L_{i-2}L_{i-1} [[a_{i-1}^+,a_i^-],[a_{i-1}^+,a_{i-
2}^-]]_q
 \hfill\break
=q\L_{i-2}L_{i-1} ([[[a_{i-1}^+,a_i^-],a_{i-1}^+]_{q^2},a_{i-2}^-]_\q
-q^2[a_{i-1}^+,[[a_i^-,a_{i-1}^+],a_{i-2}^-]_\q]_{\q^2}) = 0 $
from $(21d)$.

\n (iv) The case $i=j-1$ is similar to (iii).

Eq. $(11c)$, corresponding to $i=1$ or to $j=1$, can be proved in a 
similar manner. What remains to be proved next, are the Serre 
relations. 
Let us consider the first relation in Eq.$(12a)$.

\n (i) $i \ne 1$. 
$[e_i,e_j]=[[a_{i-1}^+,a_i^-]L_{i-1},[a_{j-1}^+,a_j^-]L_{j-1}] $
\hfill\break  
$=q^{\delta_{i-1,j}-\delta_{i-1,j-1}}
[[a_{i-1}^+,a_i^-],[a_{j-1}^+,a_j^-]]_{q^{\delta_{i,j-1}-\delta_{i-
1,j-1}}}
L_{i-1}L_{j-1} $  using $(21b)$
\hfill\break
$=([[[a_{i-1}^+,a_i^-],a_{j-1}^+]_q,a_j^-]_\q
+q[a_{j-1}^+,[[a_{i-1}^+,a_i^-],a_j^-]_\q]_\q) L_{i-1}L_{j-1} $ (from 
(23) and
the condition $|i-j| \ne 1$)
\hfill\break
$=0$ due to $(21d)$.

\n Applying the antiinvolution (14) on the first relation 
in $(12a)$, one obtaines the second relation in $(12a)$. 

\n (ii) $i = 1$. Then $[e_1,e_j]=[a_1^-,[a_{j-1}^+,a_j^-]L_{j-1}] = 
[a_1^-,[a_{j-1}^+,a_j^-]]_q L_{j-1} $
\hfill\break
$= q[[a_j^-,a_{j-1}^+],a_1^-]_\q 
L_{j-1} = 0 $  from $(21d)$ since $j>2$. 

\n We shall now prove the other Serre relations $(12b)$.

\n (i) Let $i \ne 1$. Then 
$[e_i,[e_i,e_{i+1}]_\q]_q = 
[[a_{i-1}^+,a_i^-]L_{i-1},[[a_{i-1}^+,a_i^-]L_{i-1},[a_i^+,a_{i+1}^-
]L_i]_\q]_q 
$
\hfill\break
Let us first evaluate the inner commutator . So
 
\n  $[e_i,e_{i+1}]_\q = [[a_{i-1}^+,a_i^-],[a_i^+,a_{i+1}^-]]L_iL_{i-
1} $  using $(21b)$
\hfill\break
$=([[[a_{i-1},a_i^-],a_i^+]_q,a_{i+1}^-]_\q + 
q[a_i^+,[[a_{i-1}^+,a_i^-],a_{i+1}^-]_\q]_\q) L_i L_{i-1} $   from 
(23)
\hfill\break
$=[\bar L_i a_{i-1}^+,a_{i+1}^-]_\q L_iL_{i-1} $ (due to the triple 
relations $(21d)$.

\n The full commutator reduces to 

\n $[[a_{i-1}^+,a_i^-]L_{i-1},[\bar L_i a_{i-1}^+,a_{i+1}^-]_\q 
L_iL_{i-1}]_q = -[[a_{i-1}^+,a_{i+1}^-],[a_{i-1}^+,a_i^-]]_\q L_{i-
1}^2 $
 (applying $(21a)$ and $(21b)$)\hfill\break
$= (-[[[a_{i-1}^+,a_{i+1}^-],a_{i-1}^+]_{{\q}^2},a_i^-]_q 
+{\q}^2[a_{i-1}^+,[[a_{i+1}^-,a_{i-1}^+],a_i^-]_q]_{q^2}) L_{i-1}^2 $ 
(using identity (23))  
\hfill\break
$= 0$  from $(21d)$.

\n (ii) For $i=1$,  $ [e_1,[e_1,e_2]_\q]_q = 
[a_1^-,[a_1^-,[a_1^+,a_2^-]L_1]_\q]_q $.
 
\n Let's compute the inner commutator:
 
\n $[e_1,e_2]_\q = 
[a_1^-,[a_1^+,a_2^-]L_1]_\q = [a_1^-,[a_1^+,a_2^-]]_q L_1 $  from 
(21b)
\hfill\break
$= ([[a_1^-,a_1^+],a_2^-]_q + [a_1^+,[a_1^-,a_2^-]_q]) L_1
= [[a_1^-,a_1^+],a_2^-]_q L_1 $ (due to (21e)
\hfill\break
$= {1\over (q - \q)} [L_1-\bar L_1, a_2^-]_q L_1 = a_2^- $.
\hfill\break
Therefore  $[e_1,[e_1,e_2]_\q]_q = [a_1^-,a_2^-]_q = 0 .$ 

\n The other Serre relations are poved in a similar manner. This
completes the proof.


\vskip 32pt
\noindent
{\bf 4. Concluding remarks}
\bigskip

\noindent

We have shown that apart from the Chevalley definition, the
universal enveloping algebra $U[sl(n+1)]$ of the Lie algebra
$sl(n+1)$ and also its $q-$deformed analogue, the Hopf algebra
$U_q[sl(n+1)]$, allow alternative descriptions in terms of
generators and relations. The generators are the (deformed)
creation and annihilation operators of $sl(n+1)$. In this respect the
present investigation is along the line of the results obtained
in [12, 13, 14]. In [12] the algebra $U_q[so(2n+1)]$ was quantized
via its CAGs. The latter are directly related to the quantum
statistics: in the nondeformed case the CAGs $f_1^\pm,\l,f_n^\pm,$
of $so(2n+1)$ are para-Fermi operators [4]. A concept of deformed
para-Bose operators $b_1^\pm,\l,b_n^\pm,$ was introduced in [13];
these operators provide an alternative description of the quantum
superalgebra $U_q[osp(1/2n)]$. The quantization of all Lie
superalgebras  $osp(2n+1/2m)$, namely of all (super)algebras from
the class $B$ in the Kac classification [15], via both deformed
para-Bose and deformed para-Fermi operators was carried out in
[14]. Therefore it is not surprising that the (deformed) CAGs of
$sl(n+1)$ are related to new quantum statistics, the exclusion
statistics of Haldane [5]. Clearly the present results can be
extended first of all, to all superalgebras from the class $A$,
i.e., all Lie superalgebras $sl(n/m)$ and their $q-$deformed
analogues. The CAGs of $sl(1/n)$ were studied in [1] and later on
it was shown that they describe noncanonical quantum systems with
new, quite unconventional properties [16, 17]. Again the new
statistics is an exclusion statistics [17]. It will be
interesting to extend the present approach to all simple Lie
algebras and even further to all basic Lie superalgebras.

We have not written here the explicit action of the
comultiplication $\Delta$, the counit $\varepsilon$ and the
antipode $S$ on the CAGs.  The expressions follow immediately
from the known transformations of the Chevalley generators under
$\Delta$, $\varepsilon$ and $S$ [2, 3] and the relations (13).
For instance,
$$
\Delta a_i^-
  =[[[\l[[\Delta e_1,\Delta e_2]_{\q},\Delta e_3]_{\q}
  \l]_{\q},\Delta e_{i-1}]_{\q},\Delta e_i]_{\q}, \eqno(24)
$$
where  $\Delta e_i=e_i \otimes 1 + \k_i \otimes e_i$.
Then from (19) one obtains ($i\ne 1$)
$$
\Delta e_1=a_1^- \otimes 1 + \L_1 \otimes a_1^-,\quad
\Delta e_i= [a_{i-1}^+,a_i^-]L_{i-1} \otimes 1 + 
  L_i\L_{i-1} \otimes [a_{i-1}^+ a_i^-]L_{i-1}. \eqno(25)
$$
Inserting Eqs.(25) in (24) one obtains an "explicit"
expresion for $\Delta a_i^-$ via the CAGs.
This expression  is however very involved for large $i$. 
Moreover, it is quite assymetrical for different $a_i^\pm$. 
For instance,
$$
\Delta a_1^-=a_1^-\otimes 1 + \L_1\otimes a_1^-, \eqno(26)
$$
whereas
$$
\Delta a_2^-=a_2^-\otimes 1 + \L_2\otimes a_2^- + 
(q-\q)[a_1^+,a_2^-]\otimes a_1^- . \eqno(27)
$$
Therefore, replacing in (26) the index 1 with 2, one does not
obtain $\Delta a_2^-$. One may try to define new, simpler and
more symmetric expressions for the action of $\Delta$,
$\varepsilon$ and $S$ on the CAGs, using the available
multiparameter deformations of $U_q[sl(n+1)]$ in the coalgebra
sector [18] and fixing appropriately some of the parameters.
This is the first open problem, which we would like to state.

The second problem is to construct the Fock representations of
the relations (21), namely of the deformed CAGs of
$U_q[sl(n+1)]$. This will lead directly to new solutions for the
$g$-on statistics of Karabali and Nair [9]. To this end, one has
to find as a first step, an analogue of the triple relations (10)
in the deformed case. This is equivalent to writing down the
"commutation relations" between (almost) all Cartan-Weyl
generators, expressed via the CAGs (the expressions of the
Cartan-Weyl generators via the Chevalley generators and the
commutation relations they satisfy are known [19]). As a second
step, using the Poincar\'e-Birkhoff-Witt theorem (following from
the triple relations), one can construct the Fock representations
of $U_q[sl(n+1)]$.

\bigskip\bigskip
\n {\bf Acknowledgements}

\bigskip
\n T. Palev is grateful to Prof. C. Reina for kind
hospitality at SISSA, and thanks Dr. N. I. Stoilova for
constructive discussions. He was supported by the Grant
$\Phi-416$ of the Bulgarian Foundation for Scientific Research.

\vskip 24pt
\noindent
{\bf References}

\vskip 12pt
\settabs \+ [11] & I. Patera, T. D. Palev, Theoretical interpretation 
of the 
   experiments on the elastic \cr 
   
\+ 1. & Palev T. D. 1980 {\it Journ. Math. Phys.} {\bf 21}, 1293. \cr 

\+ 2. & Drinfeld V.G. 1986 {\it Quantum Groups}, Proc. Int. Congr. 
      Math. vol. 1, 798, Berkeley.\cr

\+ 3. & Jimbo M. 1986 {\it Lett. Math. Phys.} {\bf 11}, 247. \cr

\+ 4. & Palev T. D. 1976 Thesis, Institute of Nuclear Research and
        Nuclear Energy, Sofia; \cr
\+    & 1977 Preprint JINR E17-10550; 1979 {\it Czech. J. Phys B} 
        {\bf 29}, 91; \cr 
\+    & 1980 {\it Rep. Math. Phys.}{\bf 18}, 117, 129. \cr

\+ 5. & Haldane F.D. 1991 {\it Phys. Rev. Lett.} {\bf 67}, 937.\cr

\+ 6. & Wu Y-Sh. 1994 {\it Phys. Rev. Lett.} {\bf 73}, 922. \cr

\+ 7. & Polychronakos A.P. 1996 {\it Phys. Lett. B} {\bf 365}, 202. 
\cr
 
\+ 8. & Ilinskaya A.V., Ilinsky K.N. and  Gunn J.M.F. 1996
        {\it Nucl. Phys. B} {\bf 458}, 562. \cr
        
\+ 9. & Karabali D. and Nair V. P. 1995 {\it Nucl. Phys.} {\bf B 438},
      551. \cr 

\+ 10. & Green H. S. 1993 {\it Phys. Rev.} {\bf 90}, 270. \cr

\+ 11. & Ohnuki Y. and Kamefuchi S. 1982 {\it Quantum Field
         Theory and Parastatistics} \cr
\+     & (Berlin: Springer).  \cr

\+12.  & Palev T. D. 1994 {\it Lett. Math. Phys. 31} 151. \cr

\+13.  & Palev T. D. 1993 {\it J. Phys. A: Math. Gen.} {\bf 26} 
         L1111; \cr
\+     & Hadjiivanov L. K. 1993 {\it J. Math. Phys.} {\bf 34} 5476;\cr
\+     & Palev T. D. and Van der Jeugt J. 1995 {\it J. Phys. A:
         Math. Gen.} {\bf 28} 2605. \cr

\+14.  & Palev T. D. 1996 {\it J. Phys. A: Math. Gen.} {\bf 29} 
L171;\cr
\+     & 1996   A description of the  quantum superalgebra  
         $U_q[osp(2n+1/2m)]$ via Green \cr 
\+     & generators,{\it Preprints} IC/96/127 and SISSA-119/96/FM, 
         q-alg/9607030. \cr 

\+15.  & Kac V. G. 1978 {\it Lect.Notes Math.} {\bf 626} 597.  \cr    
 
\+16.  & Palev T. D. 1982 {\it Journ. Math. Phys.} {\bf 23} 1778. \cr
\+     & Palev T. D. and Stoilova N. I. 1994 {\it J. Phys. A:
         Math. Gen.} {27} 977 and 7387. \cr

\+17.  & Palev T D and Stoilova N I 1996  Many-body Wigner quantum 
systems\cr
\+     & {\it Preprint}  IC/96/82 (1996) and hep-th/9606011. \cr

\+18.  & Dobrev V. K. and Preeti Parashar 1993 {\it J. Phys. A:
         Math. Gen.} {26} 6991.  \cr

\+19.  & Hiroyuki Y. 1989 {\it Publ. RIMS Kyout Univ.} {\bf 25} 
         503.\cr

\end